\documentclass[a4paper,11pt]{article}
\usepackage{mathrsfs}
\usepackage{latexsym}
\usepackage{amsmath,amssymb}
\usepackage{amsthm}
\usepackage{amsfonts}
\usepackage[usenames]{color}
\usepackage{amssymb}
\usepackage{graphicx}
\usepackage{amsmath}
\usepackage{amsfonts}
\usepackage{amsthm}
\usepackage{mathrsfs}
\usepackage{dsfont}
\usepackage{indentfirst}
\usepackage{enumerate}
\usepackage[square, comma, sort&compress, numbers]{natbib}

\def\lan{\langle}
\def\ran{\rangle}
\def\p{\partial }

\def\phi{\varphi}

\newtheoremstyle{mythm}{1.5ex plus 1ex minus .2ex}{1.5ex plus 1ex
minus .2ex}{\kai}{\parindent}{\song\bfseries}{}{1em}{}
\numberwithin{equation}{section}
\newtheorem{definition}{Definition}[section]
\newtheorem{theorem}{Theorem}[section]
\newtheorem{lemma}{Lemma}[section]

\newtheorem{proposition}{Proposition}[section]
\newtheorem{remark}{Remark} [section]
\newtheorem{corollary}{Corollary}[section]
\allowdisplaybreaks[4]

\begin{document}
\title
{Uniqueness of positive solutions for finsler $p$-Laplacian equations with polynomial non-linearity}
\author{Rongxun He  and  Wei Ke\footnote{Corresponding author. \textit{Email address:} wke25@swun.edu.cn.}}
\date{}
\maketitle

\begin{abstract}
 We consider the uniqueness of positive solutions to the following anisotropic elliptic equation:
\begin{equation}\nonumber
	\left\{
	\begin{aligned}
		-\Delta^F _p u&=u^q \quad \text{in} \quad  \Omega,\\
		u&=0 \quad  \text{on} \quad  \partial \Omega,  \\
	\end{aligned}
	\right.
\end{equation}
where $p>\frac{3}{2}$ and $q>p-1$. We utilize the linearized method to derive the uniqueness results, which extend the conclusion obtained by L. Brasco and E. Lindgren.
\end{abstract}

\section{Introduction}
 In this paper, we study the uniqueness of positive solutions of finsler $p$-Laplacian Dirichlet problem
\begin{equation}\label{eq2}
\left\{
\begin{aligned}
-\Delta ^F_p u&=u^q \quad \text{in} \quad  \Omega,\\
       	u&=0 \quad  \text{on} \quad  \partial \Omega,  \\
\end{aligned}
\right.
\end{equation}
where $\Omega$ is a bounded $C^2$ domain in $\mathbb{R}^N$. Finsler $p$-Laplacian operator is defined by
\begin{align*}
\Delta ^F_p u:=\text{div}(F^{p-1}(\nabla u)\nabla _{\xi}F(\nabla u)),
\end{align*}
where $F$ is a one-homogeneous positive convex $C^{3,\alpha}_{loc}(\mathbb{R}^N\setminus\{0\})$ function with $\alpha\in(0,1)$. We assume the following set of hypotheses on $F$:

\noindent \textbf{Assumption F1.} There exist $0< a_0\le b_0<+\infty$ such that  $a_0|\xi|\le F(\xi)\le b_0|\xi|$ for any $\xi\in \mathbb{R}^N$.

\noindent \textbf{Assumption F2.} There exist $0< c_0\le d_0<+\infty$ such that $c_0|\xi|^{p-2}I_N \le\nabla ^2F^p(\xi)\le d_0|\xi|^{p-2}I_N$, where $I_N$ is the identity matrix of order $N$.

\noindent \textbf{Assumption F3.} $F(t\xi)=|t|F(\xi)$ holds for any $\xi\in \mathbb{R}^N, t\in\mathbb{R}$.

We always set $H=F^p$ and denote $S^{N-1}$ the unit sphere in $\mathbb{R}^N$.

\vspace{0.5em}
\noindent \textbf{1.1. Previous results.}
The uniqueness of positive solutions to elliptic equations is an important problem in elliptic PDEs theory and has been extensively studied in the past decades. In the following, we briefly review some existing results.

We start with the Lane-Emden equation
\begin{equation}
	\label{eq0001}
	\left\{
	\begin{aligned}
		-\Delta u&=u^p \quad \text{in}  \quad \Omega,\\
		u&=0 \quad \text{on} \quad \partial \Omega.
	\end{aligned}
	\right.
\end{equation}
For the sublinear case $0<p<1$, Brezis and Oswald deduced the uniqueness of positive solutions to equation $(\ref{eq0001})$ in \cite{brezis1986uniqueness}. We also refer to a recent work \cite{BPZ2022uniqueness} for the same result without any regularity assumption on $\Omega$. For the superlinear case $1<p<\frac{N+2}{N-2}$, there are more complicated results. When $\Omega$ is a ball, Gidas, Ni and Nirenberg inferred that the problem $(\ref{eq0001})$ has exactly one radially symmetric solution in \cite{gidas1979radial}. The uniqueness result holds true if $\Omega$ is close to a ball in the sense of Hausdorff distance, see \cite{zou1994uniqueness}. For dimension $N=2$, Lin derived the uniqueness of least energy solutions to $(\ref{eq0001})$ for planar convex domain in \cite{lin1994uniqueness} and the least energy assumption can be removed if $p$ is sufficiently large, see \cite{de2019morse}. However, there are also examples that uniqueness is not correct. In \cite{Li1990multiple}, Li showed the existence of many nonradial positive solutions in an annulus of $\mathbb{R}^N$. Another fundamental counter-example is given by Dancer in \cite{Dancer1988multiple}, where $\Omega$ is an approximation to a finite union of disjoint balls.

When turning to the quasilinear Lane-Emden equation
\begin{align}
-\Delta_pu=u^{q-1} \quad\text{in}\quad\Omega,\quad\text{for}\quad 1<p<\infty,
\label{eq0002}
\end{align}
uniqueness of positive solutions holds true again for the sublinear case $1<q<p$ by \cite{diaz1987existence} (see also \cite{BPZ2022uniqueness} for general open domains). For the superlinear case $p<q<p^*$, uniqueness result holds in balls as before (see \cite{yadava1994elementary}), but fails in general. We can refer to \cite{Nazarov2000multiple} for the counter-example, which is the quasilinear version of the result in \cite{Li1990multiple}. Recently in \cite{brasco2023uniqueness}, Brasco and Lindgren proved that equation \eqref{eq0002} admits a unique ground state solution when $p>2$ and $p<q<q_0$ for some constant $q_0=q_0(N,\Omega,p)>p$. In \cite{Ke2023uniqueness}, Ke used their method to derive the uniqueness results of positive solutions to the following equation with homogeneous boundary data,
\begin{align*}
-\Delta_p u=\lambda u^{p-1}+u^{q-1}\quad\text{in}\quad\Omega,
\end{align*}
where $\lambda\in(0,\lambda_1(\Omega))$, $p>2$ and $p<q<q_0$ for some constant $q_0=q_0(N,\Omega,p,\lambda)>p$. Here, $\lambda_1(\Omega)$ denotes the first eigenvalue of $-\Delta_p$ on $\Omega$.

\vspace{0.5em}
\noindent \textbf{1.2. Main results.}
In this article, we will study the uniqueness of positive solutions to finsler $p$-Laplacian equation \eqref{eq2}. Our results are as follows.
\begin{theorem}\label{thm1}
	Let $\frac{3}{2}<p<2$, $N\geq 2$ and $p-1<q<Q_{0,p}$. Let $\Omega\subset\mathbb{R}^N$ be an open bounded connected domain of class $C^2$ such that the Dirichlet problem
	\begin{equation}\label{eq04}
		\left\{
		\begin{aligned}
			-\Delta_p  u&=u^{q} \quad \text{in} \quad  \Omega\\
			u&=0 \quad  \text{on} \quad  \partial \Omega  \\
		\end{aligned}
		\right.
	\end{equation}
has a unique positive solution $u$ that is non-degenerate.

Then, there exists a constant $\delta_1 =\delta_1 (N,\Omega,p,q)>0$ such that for every $F$ satisfying $||F(\xi)-|\xi|||_{C^{3,\alpha}(S^{N-1})}\le \delta_1$, the problem $(\ref{eq2})$ has a unique positive solution.
\end{theorem}

\begin{theorem}\label{thm2}
	Let $\frac{3}{2}<p<2$, $N\geq 2$ and let $\Omega\subset\mathbb{R}^N$ be an open bounded connected domain of class $C^2$. Then, there exists a constant $q_0=q_0(N,\Omega,p)>p-1$ such that the problem $(\ref{eq2})$ has a unique positive solution for every $q\in (p-1,q_0)$.
\end{theorem}

Theorem \ref{thm2} is an extension of the result by Brasco and Lindgren, which shows the uniqueness of least energy solutions to \eqref{eq0002} for the case $p>2$. Indeed, Theorem \ref{thm1} and \ref{thm2} still hold true for $p>2$; one can repeat all basic processes in \cite{brasco2023uniqueness} but apply the blow-up argument in \cite[Lemma 3]{lin1994uniqueness} to deal with general positive solutions. In our work, we only consider the case $p<2$ for simplicity. We remark that the corresponding problem remains open for the case $1<p\le\frac32$.

\vspace{0.5em}
\noindent  \textbf{1.3. Comments on the proof.}
The proof of Theorem \ref{thm2} is largerly inspired by that of \cite[Theorem 1.1]{brasco2023uniqueness}, dealing with the case of $p$-Laplacian for $p>2$. We utilize the linearization argument and discuss the weighted eigenvalue problem as in \cite{brasco2023uniqueness}. In the following we sketch the idea of the proof and refer to \cite{brasco2023uniqueness, lin1994uniqueness} for its root.

We argue by contradiction and assume that there exists a sequence $\{q_k\}_{k\in \mathbb{N}}$ tending to $p-1$ such that the equation $(\ref{eq2})$ admits two linearly independent positive solutions $u_k$ and $v_k$. Define $\tilde{u}_k=\frac{u_k}{M_k}$, $\tilde{v}_k=\frac{v_k}{M_k}$, where $M_k:=||u_k||_{L^{\infty}(\Omega)}\ge||v_k||_{L^{\infty}(\Omega)}$. Using the Picone-type inequality, the difference $\tilde{u}_k-\tilde{v}_k$ must be sign-changing in $\Omega$. Furthermore, it solves the linearized equation
\begin{align}
-\text{div}(A_k\nabla (\tilde{u}_k-\tilde{v}_k))=M^{q_k+1-p}_kq_k\int_{0}^{1}\big[t\tilde{u}_k+(1-t)\tilde{v}_k\big]^{q_k-1}(\tilde{u}_k-\tilde{v}_k)dt,
\label{eq0003}
\end{align}
where 
\[A_k=\frac{1}{p}\int_{0}^{1}D^2F^p(t\nabla\tilde{u}_k+(1-t)\nabla\tilde{v}_k)dt.
\]
Noticing that $A_k$ satisfies
\[\frac1C(|\nabla\tilde{u}_k|+|\nabla\tilde{v}_k|)^{p-2}|\xi|^2 \le \lan A_k\xi,\xi\ran \le C(|\nabla\tilde{u}_k|+|\nabla\tilde{v}_k|)^{p-2}|\xi|^2
\]
for every $\xi\in\mathbb{R}^N$, hence $A_k$ is singular elliptic due to $p<2$. Then we denote by $w_k$ the normalization of $\tilde{u}_k-\tilde{v}_k$ with its weighted $L^2$-norm (see \eqref{eenew111}), which infers that
\begin{align*}
\int_{\Omega}(|\nabla\tilde{u}_k|+|\nabla\tilde{v}_k|)^{p-2}|\nabla w_k|^2dx<+\infty.
\end{align*}
This enlightens us to establish a suitable compact embedding for weighted Sobolev spaces of functions satisfying
\begin{align*}
\int_{\Omega}(|\nabla\tilde{u}_k|+|\nabla\tilde{v}_k|)^{p-2}|\nabla \varphi|^2dx<+\infty.
\end{align*}
We note that in \cite{brasco2023uniqueness}, this embedding for $p>2$ plays an essential role in the proof of the main result and it takes complex computation, see \cite[Section 2.4 \& Appendix C]{brasco2023uniqueness}. However, in our case $p<2$, the method is totally different. We apply the Hopf's lemma and the Hardy's inequality to deal with $L^\sigma$-norm with weights singular on the boundary. Therefore, we can derive a uniform compact embedding into some weighted $L^\sigma$ space. 

Next, we need to explain why we consider the solutions after division. As stated in \cite{brasco2023uniqueness}, it is difficult to build a uniform estimate for solutions directly. Instead, we let the maximum value of solutions divide the equation and thus can apply the uniform $C^{1,\alpha}$ estimates by classical results of degenerate elliptic equations. To guarantee the constant $M_k^{q_k+1-p}$ in \eqref{eq0003} is bounded uniformly for $\{q_k\}_{k\in\mathbb{N}}$, we learn from \cite[Lemma 3]{lin1994uniqueness} and make use of the blow-up method and Liouville's Theorem. Moreover, we can prove $M_k^{q_k+1-p}$ converges to the first eigenvalue of $-\Delta_p^F$ on $\Omega$ as $q_k\to p-1$.

Finally, we pass the limit in the linearized equation for $w_k$ and obtain a convergence to a non-trivial limit function $w$ thanks to compact embeddings. Furthermore, we can prove $w$ belongs to $X_0^{1,2}(\Omega;|\nabla\tilde u|^{p-2})$, which denotes the completion of $C_c^\infty(\Omega)$ with respect to the norm
\begin{align*}
\Vert\varphi\Vert_{X^{1,2}(\Omega;|\nabla\tilde u|^{p-2})}=||\phi||_{L^2(\Omega)}+\left(\int_\Omega|\nabla\tilde u|^{p-2}|\nabla\varphi|^2dx\right)^{\frac12}.
\end{align*}
Here $\tilde u$ denotes the limit function of $\tilde u_k$ as $k\to\infty$. To ensure every $\varphi\in C_c^\infty(\Omega)$ satisfies $\Vert\varphi\Vert_{X^{1,2}(\Omega;|\nabla\tilde u|^{p-2})}<\infty$, we utilize the striking regularity results in \cite{castorina2019hopf} (see \cite{brasco2023uniqueness, damascelli2004} for the case of $p$-Laplacian), which show the integrablity of $|\nabla\tilde u|$ raised with some suitable negative powers. Indeed, for every $r<p-1$, there exists a constant $\mathcal{S}>0$ such that
\begin{align}
\int_\Omega\frac1{|\nabla\tilde u(x)|^r}dx\le \mathcal{S}.
\label{eq0004}
\end{align}
We apply \eqref{eq0004} for the choice $r=-(p-2)$. Then $r<p-1$ infers $p>\frac32$; this is the main reason of the condition $\frac{3}{2}<p<2$ in our theorems. We conclude the proof by a similar argument of \cite{brasco2023uniqueness}: obtain that $w\in X_0^{1,2}(\Omega;|\nabla u|^{p-2})$ is a sign-changing first eigenfunction of a weighted linear eigenvalue problem. This leads to a contradiction to the sign property of the first eigenfunction.

\vspace{0.5em}
\noindent \textbf{1.4. Plan of the paper.} 
In Section 2, we prove some preliminary results. In particular, we devote this section to proving uniform $L^\infty$-estimate and weighted embeddings. In Section 3, we study the first eigenvalue and first eigenfunction of a weighted linear eigenvalue problem. They are suitable adaptation of \cite[Section 3]{brasco2023uniqueness} and we decide to keep the details for the sake of completeness and for the benefit of the reader. In Section 4, we prove Theorem \ref{thm2}, along the lines as above. Finally in Section 5, we prove Theorem \ref{thm1}.

\vspace{0.5em}
\noindent\textbf{Acknowledgements.} The authors are grateful to Professor Genggeng Huang for suggesting this problem
and for helpful discussions.

\section{Preliminaries}
In this section, we will introduce some preliminary results for latter applications. It is worth pointing out that these results applies to all $p>1$.

\vspace{0.5em}
\noindent \textbf{2.1. $L^\infty$-estimate.}
In the following, we prove the uniform a priori $L^\infty$-estimate for finsler $p$-Laplacian equations. We first give a classical result of boundary gradient estimates. For the case of Laplacian, we refer to \cite[Theorem 1.1.14]{han2016nonlinear}. 
\begin{lemma}\label{lem1}
   Let $\Omega\subset \mathbb{R}^N$ be a domain satisfying the uniform exterior sphere condition; i.e., there is $r_0>0$ such that, for every $x^0\in \partial \Omega$, there is $y^0\in \mathbb{R}^N \backslash \Omega$ with $\overline{B_{r_0}(y^0)}\cap \overline{\Omega}={x^0}$. Let $u\in  C(\overline{\Omega})\cap C^{1,\alpha}(\Omega)$ be a solution of
   \begin{equation}\label{eq5}
   	\left\{
   	\begin{aligned}
   		-\Delta ^F_p u&=f \quad \text{in} \quad  \Omega\\
   		u&=0 \quad  \text{on} \quad  \partial \Omega  \\
   	\end{aligned}
   	\right.
   \end{equation}
with $f\in L^{\infty}(\Omega)\cap C(\overline{\Omega})$ and $\alpha\in(0,1)$. Let $M>0$ satisfy $||u||_{L^{\infty}(\Omega)}+||f||_{L^{\infty}(\Omega)}<M$. Then for any $x\in \Omega\cap B_1(x^0)$, $|u(x)-u(x^0)|\leq C|x-x^0|$, where $C$ is a positive constant depending only on $p$, $N$, $r_0$, $M$ and assumptions on $F$.
\end{lemma}
\noindent \emph{Proof.}  For given $x^0\in \partial \Omega$ and $y^0\in \mathbb{R}^N \backslash \Omega$, let $d(x)$ be the distance from $x$ to $\partial  B_{r_0}(y^0)$ for $x\in\Omega$, i.e., $d(x):=|x-y^{0}|-r_0$. 

Claim: there exists a $C^2$ function $\psi$ defined on $[0,1)$, with $\psi(0)=0$ , $\psi'>0$ in $(0,1)$ and such that $\omega=\psi\circ d$ satisfies
\[\Delta^F_p\omega<-C(p,N,r_0)<0\quad \text{in} \quad\Omega\cap B_1(x^0).\]
For the proof of the claim, we first compute
\begin{gather*}
\partial_i d(x)=\frac{x_i-y^0_{i}}{|x-y^0|},\\
\partial_{ij} d(x)=\frac{\delta_{ij}}{|x-y^0|}-\frac{(x_i-y^0_{i})(x_j-y^0_{j})}{|x-y^0|^3}.    
\end{gather*}
Then
\[\Delta ^F_p \omega =\text{div}(F^{p-1}(\nabla \omega)F_{\xi}(\nabla \omega))=F^{p-2}\sum_{i,j}\left((p-1)F_{\xi_i}F_{\xi_j}+FF_{\xi_i\xi_j}\right)(\psi''\partial_id\partial_jd+\psi '\partial _{ij}d).
\]
We now require $\psi''<0$.
Set $a=\frac{Nd_0-c_0}{r_0c_0}>0$, and by recalling Assumption F2, we get
\[\Delta _p^F\omega \le \frac{c_0}{p}(\psi')^{p-2}(\psi ''+a\psi').
\]
Then we choose a positive constant $b$ and find a function $\psi$ in $[0,1)$ such that
\[\psi''+a\psi'+b=0\quad \text{in} \quad (0,1).
\]
A particular solution of the previous ordinary differential equation is given by
\[\psi(d)=-\frac{b}{a}d+\frac{A}{a}\left(1-e^{-ad}\right)
\]
for some constant $A$. Set $A=\frac{2b}{a}e^{a}>0$, then
$\psi$ satisfies all the requirements we imposed.
Then we get
\begin{equation*}
	\left\{
	\begin{aligned}
		\Delta ^F_p(C\omega)&\leq\Delta ^F_p  u \quad \text{in} \quad \Omega\cap B_1(x^0),\\
		C\omega&\geq  u \quad \text{on} \quad \partial (\Omega \cap B_1(x^0)),
	\end{aligned}
\right.
\end{equation*}
where $C$ is a positive constant depending only on $p$, $N$, $r_0$, $c_0$, $d_0$ and $M$. By the comparison principle, we obtain 
\[u\leq C\omega \quad \text{in} \quad \Omega \cap B_1(x^0).
\]
Similarly we have
\[-u\leq C\omega \quad \text{in} \quad \Omega \cap B_1(x^0).
\]
Recalling that $u(x^0)=\omega(x^0)=0$, we obtain
\[|u-u(x^0)|\leq C|\omega-\omega(x^0)| \quad \text{in} \quad \Omega \cap B_1(x^0).
\]
This implies the desired result. \hfill $ \Box$
\begin{remark}
	\label{rem2}
	Let $\Omega$ satisfy the uniform exterior sphere condition with the radius $r_0$. Then for any $\mu>1$, the domain $\mu\Omega$ satisfies the uniform exterior sphere condition in particular with the same radius $r_0$.
\end{remark}

Next, we will derive the a priori $L^\infty$-estimate for solutions to $(\ref{eq2})$. Here we make use of the blow-up argument by Lin, see \cite[Lemma 3]{lin1994uniqueness}.
\begin{lemma}
	\label{lem2} 
	Let $\Omega\subset \mathbb{R}^N$ be an open bounded connected set of class $C^2$. Let $p>1$ and $u_q\in W^{1,p}_0(\Omega)$ be a positive solution of $(\ref{eq2})$. Set $M_q=||u_q||_{L^{\infty}(\Omega)}$. Then, there exists a $\delta_0=\delta_0(p,N,\Omega)>0$, such that $M_q^{q+1-p}\leq C$ for $ q\in [p-1,p-1+\delta_0]$, where the constant $C$ depends only on $p,N,\delta_0$ and $\Omega$.
\end{lemma}
\noindent \emph{Proof.} Suppose on the contrary that there exists $\{q_k\}\rightarrow p-1$, such that $M_{q_k}^{q_k+1-p}=||u_k||_{L^{\infty}(\Omega)}^{q_k+1-p}\rightarrow \infty$, where $u_k$ is a positive solution of $(\ref{eq2})$ with $q=q_k$. Let $\{x^k\}_{k\in \mathbb{N}}$ be a sequence in $\Omega$ such that $u_k(x^k)=M_{q_k}$. Define
\[\tilde{u}_k(y)=\frac{u_k(\mu_ky+x^k)}{M_{q_k}}.\]
Then $\tilde{u}_k$ solves
\begin{equation}
	\label{eq13}
    \left\{
    \begin{aligned}
    	-\Delta ^F_p \tilde{u}_k&=M_{q_k}^{q_k+1-p}\mu^p_k\tilde{u}^{q_k}_k  \quad \text{in} \quad \Omega_k,\\
    	\tilde{u}_k&=0 \quad \text{on} \quad \partial \Omega_k,
    \end{aligned}
    \right.
\end{equation}
where $\Omega_k=\{y\in \mathbb{R}^N:\mu_ky+x^k\in \Omega\}. $ We have $0\leq\tilde{u}_k\leq1$, $\tilde{u}_k(0)=1$.
Let $M_{q_k}^{q_k+1-p}\mu^p_k=1$, then $\mu_k\rightarrow0$ as $ k\rightarrow\infty$. Equation \eqref{eq13} turns into
\begin{equation*}
	\left\{
	\begin{aligned}
		-\Delta^F_p \tilde{u}_k&=\tilde{u}^{q_k}_k  \quad \text{in} \quad \Omega_k,\\
		\tilde{u}_k&=0 \quad \text{on} \quad \partial \Omega_k.
	\end{aligned}
	\right.
\end{equation*}
Up to a subsequence, two situations may occur:
\[\text{either} \quad \text{dist}(x^k,\partial \Omega)\mu^{-1}_k\rightarrow +\infty \quad \text{or} \quad \text{dist}(x^k,\partial \Omega)\mu^{-1}_k\rightarrow d \geq 0.
\]

\textbf{Case 1:} $\lim_{k\rightarrow +\infty}\text{dist}(x^k,\partial \Omega)\mu^{-1}_k\rightarrow +\infty$. Then for any $r>1$, $B_r(0)\subset\Omega_k$ holds for $k$ large. According to the local $C^{1,\alpha}$ estimates (see \cite[Theorem 1 \& 2]{dibenedetto1983c1+}), there exists a constant $C=C(r)>0$ uniformly for $k$ such that 
\[||\tilde{u}_k||_{C^{1,\alpha}(B_r(0))}\le C.
\]
By selecting a subsequence, $\tilde{u}_k\rightarrow \tilde{u}$ in $C^{1,\beta}_{loc}(\mathbb{R}^N)$ for some $\beta\in (0,\alpha)$.
Moreover, $\tilde{u}(0)=1$, $\tilde{u}\ge 0$ in $\mathbb{R}^N$ and $\tilde{u}$ solves
\[\int_{\mathbb{R}^N}F^{p-1}(\nabla \tilde{u})\nabla_{\xi} F(\nabla \tilde{u})\cdot \nabla \psi dx=\int_{\mathbb{R}^N}\tilde{u}^{p-1}\psi dx \quad \text{for all} \quad \psi \in C^\infty_c(\mathbb{R}^N),
\]
namely, $\tilde{u}$ is a distributional solution of
\begin{equation}
	\label{eq17}
	-\Delta^F_p \tilde{u}=\tilde{u}^{p-1} \quad \text{in} \quad \mathbb{R}^N.
\end{equation}
We claim that $\tilde{u}=0$. Let $\varphi _R$ be the first eigenfunction of $-\Delta^F_p$ on $B_R(0)$, i.e., $\varphi_R$ satisfies 
\begin{equation}
	\label{eq16}
	\left\{
	\begin{aligned}
		-\Delta^F_p \varphi_R&=\lambda_R \varphi^{p-1}_R \quad \text{in} \quad B_R(0),\\
		\varphi_R&=0 \quad \text{on} \quad \partial B_R(0),
	\end{aligned}
	\right.
\end{equation}
where $\lambda_R$ is the first eigenvalue of $-\Delta^F_p$ on $B_R(0)$ and $\lambda_R\rightarrow 0$ as $R\rightarrow +\infty$. If $\tilde u\neq0$, then by testing $\psi_R=\varphi_R^p\tilde{u}^{1-p}$ on $(\ref{eq17})$, we obtain
\[\int_{B_R(0)}F^{p-1}(\nabla \tilde{u})\nabla_{\xi}F(\nabla \tilde{u})\cdot \nabla(\varphi_R^p\tilde{u}^{1-p})dx= \int_{B_R(0)}\varphi^p_R dx.
\]
Denote the left-hand-side by LHS. Then, by Young's inequality,
\begin{equation*}
	\begin{split}
		\text{LHS}=&(1-p)\int_{B_R(0)}\varphi_R^pF^p(\nabla \tilde{u}) \tilde{u} ^{-p}dx+p\int_{B_R(0)}F^{p-1}(\nabla \tilde{u})\nabla _{\xi}F(\nabla \tilde{u})\tilde{u}^{1-p}\varphi^{p-1}_R\nabla \varphi_Rdx\\
		\le& (1-p)\int_{B_R(0)}\varphi_R^pF^p(\nabla \tilde{u}) \tilde{u} ^{-p}dx +(p-1)\int_{B_R(0)}\varphi_R^pF^p(\nabla \tilde{u}) \tilde{u} ^{-p}dx\\
		&+C_p\int_{B_R(0)}|\nabla \varphi_R|^pdx.
	\end{split}
\end{equation*}
Thus we have
\[\int_{B_R(0)}|\nabla \varphi_R|^pdx\ge C_p^{-1}\int_{B_R(0)}|\varphi_R|^pdx.
\]
By $(\ref{eq16})$ and Assumption F1, 
\[\lambda_R\int_{B_R(0)}|\varphi_R|^pdx \ge C_p^{-1} a_0^p\int_{B_R(0)}|\varphi_R|^pdx.
\]
Let $R\rightarrow \infty$, and we reach a contradiction. Thus $\tilde{u}\equiv 0$.
But this contradicts $\tilde{u}(0)=1$.

\vspace{0.5em}
\textbf{Case 2:} $\text{dist}(x^k,\partial \Omega)\mu_k^{-1}\rightarrow d\ge 0$. Then $x^k\rightarrow x^0\in \partial \Omega$. Let $d_k$ be the distance from $x^k$ to $\partial \Omega$ and $d_k=\text{dist}(x^k,\tilde{x}^k)$, where $\tilde{x}^k\in \partial \Omega$. After a translation and rotation, we may assume $\tilde{x}^k=0$, $\nu(0)=-e_N$, where $\nu(0)$ is the outward unit normal vector at $0$. 

Claim: $d>0$. If this is not true, we define the function $v_k(y)=u_k(\mu_ky)/M_{q_k}$ for $y\in \tilde{\Omega}_k$, where  $\tilde{\Omega}_k=\{y\in \mathbb{R}^N:\mu_ky\in \Omega\}$. Then $v_k$ satisfies
\begin{equation}
	\label{eq23}
	\left\{
	\begin{aligned}
		-\Delta^F_p v_k&=v^{q_k}_k  \quad \text{in} \quad \tilde{\Omega}_k,\\
		v_k&=0 \quad \text{on} \quad \partial \tilde{\Omega}_k.
	\end{aligned}
	\right.
\end{equation}
Setting $z^k={x^k}/{\mu_k}$, we have $v_k(z^k)=1$ and 
\[z^k\rightarrow 0 \quad \text{as} \quad k\rightarrow \infty.
\]
By Lemma $\ref{lem1}$ and Remark $\ref{rem2}$, for $z^k\in \tilde{\Omega}_k\cap B_1(0)$, we have
\[|v_k(z^k)-v_k(0)|\le C|z^k|\rightarrow  0 \quad \text{as} \quad k \rightarrow \infty.
\]
But this contradicts $v_k(z^k)=1$. Hence, $d>0$ and $z^k\rightarrow de_N$ as $k\rightarrow\infty$.

Denote $x=(x',x_N)$ where $x'=(x_1,\cdots,x_{N-1})$. Then near the origin, $\partial \Omega$ can be represented as $x_N=\rho(x')$ with $\rho (0)=0$ and $\nabla \rho(0)=0$. After the scaling $y=x/\mu_k$, we have $y_N={\rho(\mu_ky')}/{\mu_k} \rightarrow 0 $ locally uniformly as  $k\rightarrow \infty$.

Using the local $C^{1,\alpha}$ estimates again, we have $v_k\rightarrow v$ in $C^{1,\beta}_{loc}(\mathbb{R}^N_{+})$, up to a subsequence. Moreover, an application of local gradient estimates yields that $v$ vanishes on the boundary $\partial\mathbb{R}^N_{+}$. Therefore, by passing the limit in \eqref{eq23}, we have
\begin{equation*}
	\left\{
	\begin{aligned}
		-\Delta^F_p v&=v^{p-1} \quad \text{in} \quad \mathbb{R}^N_+,\\
		v&=0 \quad \text{on} \quad \partial\mathbb{R}^N_{+},\\
        v(de_N)&=1.\\
	\end{aligned}
	\right.
\end{equation*}
Let $\varphi_R$ be the first positive eigenfunction of $-\Delta^F_p$ on $B_R(Re_N)$. Then a similar argument as in Case 1 infers $v\equiv 0$. But this contradicts $v(de_N)=1$. The lemma is proved. \hfill $\Box$

\begin{lemma}
	\label{lem3}
    Suppose the same conditions as in Lemma \ref{lem2}. Let $\{q_k\}\rightarrow p-1$ and $\{u_k\}$ be a sequence of positive solutions to $(\ref{eq2})$ with $q=q_k$. Then up to a subsequence, we have $u_k/{||u_k||_{L^\infty(\Omega)}} \rightarrow \tilde{u}$ in $C^{1,\beta}(\overline\Omega)$ for some $\beta\in (0,1)$, and $\tilde{u}$ is the first eigenfunction of $-\Delta^F_p$ on $\Omega$.
\end{lemma}
\noindent \emph{Proof.} Denote $M_{q_k}:=||u_k||_{L^{\infty}(\Omega)}$ and $\tilde{u}_k=u_k/M_{q_k}$. Then $\tilde{u}_k$ solves
\begin{equation*}
	\left\{
	\begin{aligned}
		-\Delta_p^F \tilde{u}_k&=M^{q_k+1-p}_{q_k}\tilde{u}^{q_k}_k \quad \text{in} \quad  \Omega,\\
		\tilde{u}_k&=0 \quad \text{on} \quad \partial \Omega.
	\end{aligned} 
	\right.
\end{equation*}
By Lemma \ref{lem2}, we have $M^{q_k+1-p}_{q_k}\leq C$. Then using the global $C^{1,\alpha}$ estimates (see \cite[Theorem 1]{lieberman1988boundary}), we have the uniform estimates $||\tilde{u}_k||_{C^{1,\alpha}(\overline\Omega)}\le C$. By selecting a subsequence, $\tilde{u}_k\rightarrow \tilde{u}\ge 0$ in $C^{1,\beta}(\overline{\Omega})$ for some $\beta\in(0,\alpha)$ and $||\tilde{u}||_{L^{\infty}(\Omega)}=1$. Moreover, $\tilde{u}$ satisfies
\begin{equation*}
	\left\{
	\begin{aligned}
		-\Delta^F_p \tilde{u}&=C_1\tilde{u}^{p-1} \quad \text{in} \quad  \Omega,\\
		\tilde{u}&=0 \quad \text{on} \quad \partial \Omega.
	\end{aligned} 
	\right.
\end{equation*}  
Since the uniqueness of the (first) eigenvalue problem of $-\Delta^F_p$, $\tilde u$ is the first eigenfunction and $C_1$ is the first eigenvalue of $-\Delta^F_p$ on $\Omega$.  \hfill $ \Box$

\vspace{0.5em}
\noindent \textbf{2.2. Uniform weighted embeddings.} We next consider the following equation
\begin{equation}
	\label{eq27}
	\left\{
	\begin{aligned}
		-\Delta^F_p u&=\mu u^{q}\quad \text{in} \quad  \Omega,\\
		u&=0 \quad \text{on} \quad \partial \Omega,
	\end{aligned}
	\right.
\end{equation}
where $\mu\in[0,\hat{\mu}]$, $\hat{\mu}>0$ and $q\ge p-1$. We will obtain the Hardy-type inequality and uniform weighted embeddings for positive solutions of \eqref{eq27}.

Before that, we first state the following regularity result. The proof is based on compactness argument and we refer to \cite[Theorem 2.5]{brasco2023uniqueness} for details.
\begin{theorem}
	\label{thm5}
     Let $\Omega\subset \mathbb{R}^N$ be an open bounded connected set of class $C^2$. Let $p>1$ and $u_{F,q}\in W^{1,p}_0(\Omega)$ be a positive solution of $(\ref{eq27})$ with $||u_{F,q}||_{L^{\infty}(\Omega)}= 1$. Then there exist positive constants $\delta_0$, $\delta_1$ such that for every $F$ with $||F(\xi)-|\xi|||_{C^{3,\alpha}(S^{N-1})}\le \delta_1$ or every $q\in[p-1,p-1+\delta_0]$, the following holds:
\begin{enumerate}[\rm(1)]
    \item $u_{F,q}\in C^{1, \gamma}(\overline{\Omega})$ with the uniform estimate
	\begin{equation}
		||u_{F,q}||_{C^{1, \gamma}(\overline{\Omega})}\le L, \nonumber
	\end{equation}
    for some $\gamma=\gamma(p, \hat\mu, \alpha, N, \delta_0, \delta_1, \Omega)\in (0,1)$, $L=L(p, \hat\mu, \alpha, N, \delta_0, \delta_1, \Omega)>0$;
    \item by defining $\Omega_{\tau}=\{x\in\Omega:\text{dist}(x,\p\Omega)\le\tau\}$, we have 
    \begin{equation}
    	|\nabla u_{F,q}|\ge \mu_0 \quad \text{in} \quad \Omega_{\tau}, \nonumber
    \end{equation}
    and
    \begin{equation}
    	 u_{F,q}\ge \mu_1 \quad \text{in} \quad \overline{\Omega \backslash \Omega_{\tau}}, \nonumber
    \end{equation}
    for some $\tau=\tau(p, \hat\mu, \alpha, N, \delta_0, \delta_1, \Omega)>0$, $\mu_i=\mu_i(p, \hat\mu,\tau, \alpha, N, \delta_0, \delta_1, \Omega)>0$ $(i= 1, 2)$.
\end{enumerate}
\end{theorem}
\begin{theorem}[Uniform Hardy-type inequality]
	\label{thm6new}
    Let $p>1$ and let $\Omega\subset \mathbb{R}^N$ be an open bounded connected set of class $C^2$. For every $F$ with $||F(\xi)-|\xi|||_{C^{3,\alpha}(S^{N-1})}\le \delta_1$ or every $q\in [p-1, p-1+\delta_0]$, let $u_{F,q}\in W^{1,p}_0(\Omega)$ be a positive solution of $(\ref{eq27})$ with $||u_{F,q}||_{L^{\infty}(\Omega)}= 1$. Then for $s_0>-2$ fixed, there exists a constant $\sigma_0>2$ such that for every $2\le\sigma\le\sigma_0$ and $s\ge s_0$, it follows that
\begin{equation*}
\mathcal{T}\left(\int_{\Omega}|u_{F,q}|^{s}|\phi|^{\sigma}dx\right)^{\frac{2}{\sigma}}\le \int_{\Omega} |\nabla \phi|^2dx \quad \text{for every} \quad \phi\in C^1(\overline{\Omega})\cap W^{1,p}_0(\Omega),
\end{equation*}
for some constant $\mathcal{T}=\mathcal{T}(p, \hat{\mu}, \alpha, N, \delta_0, \delta_1, \sigma, s, \Omega)>0$.
\end{theorem}
\noindent \emph{Proof.} We only deal with the case $s<0$. By the classical Hardy's inequality \cite{hardy1997}, we have
\begin{equation}
\int_{\Omega}\frac{\phi^2}{d^2}dx \le C\int_{\Omega}|\nabla \phi|^2 dx \quad \text{for every} \quad \phi\in C^1(\overline{\Omega})\cap W^{1,p}_0(\Omega),
\nonumber
\end{equation}
where $d(x)=\text{dist}(x,\partial\Omega)$. By Theorem \ref{thm5}, we obtain
\begin{equation}
u_{F,q}\ge\frac{\mu_0}{2}d \quad \text{ in } \Omega_\tau,
\nonumber
\end{equation}
when $\tau>0$ is sufficiently small. Moreover, $u_{F,q}\ge \mu_1$ in $\overline{\Omega \backslash \Omega_{\tau}}$. Thus it follows
\begin{equation}
\int_{\Omega} |u_{F,q}|^{-2}|\phi|^2 dx \le C\int_{\Omega} |\nabla \phi|^2 dx,
\nonumber
\end{equation}
with a uniform constant $C>0$.

Define
\[\sigma_0=2+\frac{2(2+s_0)}{N-2}>2.\]
By direct calculation,
\begin{equation}
\frac{2(\sigma+s)}{2+s}\le2^*=\frac{2N}{N-2}  \nonumber
\end{equation}
holds for $2\le\sigma\le\sigma_0$ and $s\ge s_0$.
Then by using the H\"older's inequality and Sobolev's inequality, we obtain
\begin{align}
\left(\int_{\Omega}|u_{F,q}|^{s}|\phi|^{\sigma}dx\right)^{\frac{2}{\sigma}}&=\left(\int_{\Omega}|u_{F,q}|^{s}|\phi|^{-s}|\phi|^{\sigma+s}dx\right)^{\frac{2}{\sigma}}  \nonumber\\
&\le \left(\int_{\Omega}|u_{F,q}|^{-2}|\phi|^{2}dx\right)^{\frac{-s}{\sigma}}\left(\int_{\Omega}|\phi|^{\frac{2(\sigma+s)}{2+s}}dx\right)^{\frac{2+s}{\sigma}}  \nonumber\\
&\le \left(C\int_{\Omega}|\nabla \phi|^2 dx\right)^{\frac{-s}{\sigma}}\left(C\int_{\Omega}|\nabla \phi|^{2}dx\right)^{\frac{\sigma+s}{\sigma}}  \nonumber\\
&= \tilde{C}\int_{\Omega}|\nabla \phi|^{2}dx.  \nonumber
\end{align}
We complete the proof by taking $\mathcal{T}=\tilde{C}^{-1}$.  \hfill $\Box$

\begin{remark}
\label{rmknew1}
	If $1<p<2$, then combining the above theorem and global $C^{1,\alpha}$ estimates of solutions, we have
\begin{equation}
	\label{eq34.1}
	\mathcal{T}\int_{\Omega}|u_{F,q}|^{p-2}|\phi|^{2}dx \le \int_{\Omega} |\nabla u_{F,q}|^{p-2}|\nabla \phi|^2dx
\end{equation}
holds uniformly for every $\phi\in C^1(\overline{\Omega})\cap W^{1,p}_0(\Omega)$.
\end{remark}

Using the previous Hardy-type inequality, we establish the uniform weighted compact embeddings as follows.
\begin{corollary}[Uniform weighted compact embedding]
	\label{cor2.2}
	Let $p>1$ and let $\Omega \subset \mathbb{R}^N$ be an open bounded connected set of class $C^2$. For $\{F_k\}_{{}k\in \mathbb{N}}$ with $||F_k(\xi)-|\xi|||_{C^{3,\alpha}(S^{N-1})}\le \delta_1$ or $\{q_k\}_{{}k\in \mathbb{N}}\subset [p-1,p-1+\delta_0]$, consider accordingly $u_k\in W^{1,p}_0(\Omega)$  positive solutions of $(\ref{eq27})$ with $||u_k||_{L^{\infty}(\Omega)}=1$. If $\{\phi _k\}_{k\in \mathbb{N}}\subset C^1(\overline{\Omega})\cap W^{1,p}_0(\Omega)$ satisfies
	\[\int_{\Omega}|\nabla \phi_k|^2dx\le C\quad \text{for  every} \quad k \in \mathbb{N},
	\]
	then $\{\phi_k\}_{k\in \mathbb{N}}$ converges to $\phi$ strongly in $L^2(\Omega)$ and weakly in $W^{1,2}(\Omega)$, up to a subsequence. Moreover, it follows
     \[\int_{\Omega}|u_k|^{q_k-1}|\phi_k-\phi|^2 dx\rightarrow 0 \quad \text{ as } k\rightarrow\infty.\]
\end{corollary}
\noindent \emph{Proof.} The first result is an easy consequence of the classical Rellich-Kondra\v{s}ov Theorem. For the remainder, we compute
\begin{align}
\int_{\Omega}|u_k|^{q_k-1}|\phi_k-\phi_l|^2 dx &\le \left(\int_{\Omega}|u_k|^{2(q_k-1)}|\phi_k-\phi_l|^2 dx\right)^{\frac12}\left(\int_{\Omega}|\phi_k-\phi_l|^2 dx\right)^{\frac12}   \nonumber  \\
&\le \left(\mathcal{T}^{-1}\int_{\Omega}|\nabla(\phi_k-\phi_l)|^2 dx\right)^{\frac12}\left(\int_{\Omega}|\phi_k-\phi_l|^2 dx\right)^{\frac12}   \nonumber  \\
&\le C||\phi_k-\phi_l||_{L^2(\Omega)}, \nonumber 
\end{align} 
where the second inequality follows by Theorem \ref{thm6new}. We conclude by letting $l\rightarrow\infty$.    \hfill $\Box$

\vspace{0.5em}
\noindent \textbf{2.3. Other Properties.}
We will introduce two key results of positive solutions to \eqref{eq27}. The first lemma states that the difference of any two different solutions must change sign. We refer to \cite[Lemma 2.2]{brasco2023uniqueness} for the proof; note that the general Picone's inequality is used for the case of finsler $p$-Laplacian, see for example \cite[Proposition 2.9]{brasco2014convexity}.
\begin{lemma}
	\label{lem 3.8}
	Let $p>1$ and $q>p-1$. Let $\Omega\subset \mathbb{R}^N$ be an open bounded connected set.  Let $u, v\in W^{1,p}_0(\Omega)$ be two distinct positive solutions of equation $(\ref{eq27})$. Then we must have 
	\[|\{x\in\Omega:u(x)>v(x)\}|>0 \quad \text{and} \quad |\{x\in\Omega:u(x)<v(x)\}|>0. 
	\]
	In other words, the difference $u-v$ must change sign in $\Omega$.
\end{lemma}

The second result describes the crucial regularity estimates for positive solutions of finsler $p$-Laplacian equations, namely the (weighted) integrability of the inverse of the gradient. This theorem was proved by \cite{damascelli2004} for $p$-Laplacian equations, and by \cite{castorina2019hopf} for general quasilinear anisotropic equations. We also refer to \cite{brasco2023uniqueness} for precise statements on how the estimates depend on the solution.

\begin{theorem}
\label{thmnew}
Let $p>1$ and let $\Omega \subset \mathbb{R}^N$ be an open bounded connected set of class $C^2$. For every $F$ with $||F(\xi)-|\xi|||_{C^{3,\alpha}(S^{N-1})}\le \delta_1$ or every $q\in [p-1, p-1+\delta_0]$, let $u\in W^{1,p}_0(\Omega)$ be a positive solution of equation \eqref{eq27} with $||u||_{L^\infty(\Omega)}=1$. Then for
\begin{equation}
	\left\{
	\begin{aligned}
		\gamma <N-2, \quad \text{if} \quad  N\ge 3,\\
		\gamma \le 0, \quad \text{if} \quad N=2. \nonumber
	\end{aligned}
	\right.
\end{equation}
and every $r<p-1$, there exists $\tilde{\mathcal{S}} =\tilde{\mathcal{S}}(p, \hat\mu, \alpha, N, \delta_0, \delta_1, \gamma, r, \Omega)>0$ such that 
\begin{equation}
    \label{eqnew1212}
	\sup_{x\in \Omega}\int_{\Omega}\frac{1}{|\nabla u(y)|^{r}|x-y|^{\gamma}}dy\le \tilde{\mathcal{S}}.
\end{equation}
In particular, we have
\begin{equation}
	\label{eqnew2121}
	\int_{\Omega}\frac{1}{|\nabla u(y)|^{r}}dy\leq\mathcal{S},
\end{equation}
for some $\mathcal{S}=\mathcal{S}(p, \hat\mu, \alpha, N, \delta_0, \delta_1, r, \Omega)>0$.
\end{theorem}
We remark that the estimate \eqref{eqnew1212} plays a crucial role in proving the weighted embedding \eqref{eq34.1} for the case $p>2$, see \cite{brasco2023uniqueness,damascelli2004}. However, in our case $p<2$, we only require the estimate \eqref{eqnew2121}.

\section{Eigenvalue problem}
In this section, we study a weighted eigenvalue problem for finsler $p$-Laplacian equation. As stated, this part is largely inspired by \cite[Section 3]{brasco2023uniqueness}, where the authors treat the case of $p$-Laplacian with $p>2$.

We first include the following elementary inequality.
\begin{lemma}
	\label{lem4}
	Let $p>1$ and let $\Omega\subset \mathbb{R}^N$ be an open set. Then for every $v, w, \phi\in W^{1,1}_{loc}(\Omega)$, we have
	\begin{align*}
	\big|\lan D^2H(\nabla \phi)\nabla v,\nabla v\ran&-\lan D^2H(\nabla \phi)\nabla w,\nabla w\ran\big|\\
	&\le d_0|\nabla \phi|^{p-2}|\nabla v-\nabla w|(|\nabla v|+|\nabla w|),\quad a.e. \quad \text{in} \quad \Omega.
	\end{align*}
\end{lemma}
\noindent \emph{Proof.} Observe that $D^2H$ is symmetric, then by direct computation,
\begin{align*}
	&\big|\lan D^2H(\nabla \phi)\nabla v,\nabla v\ran-\lan D^2H(\nabla \phi)\nabla w,\nabla w\ran\big| \\
    =&\big|\lan D^2H(\nabla \phi)(\nabla v-\nabla w),\nabla v\ran +\lan D^2H(\nabla \phi)(\nabla v-\nabla w),\nabla w\ran\big| \\
    \le&\big|D^2H(\nabla \phi)(\nabla v-\nabla w)\big|(|\nabla v|+|\nabla w|).
\end{align*}
By Assumption F2, we have
\begin{equation}
	\nonumber
	\big|D^2H(\nabla \phi)(\nabla v-\nabla w)\big|\le d_0|\nabla \phi|^{p-2}|\nabla v-\nabla w|.
\end{equation}
We complete the proof by inserting this inequality into the previous estimate.\hfill $\Box$

\begin{proposition}
	\label{pro3.2}
	Let $p>1$ and let $\Omega \subset \mathbb{R}^N$ be an open bounded connected set of class $C^2$. Let $U\in W^{1,p}_0(\Omega)$ be the unique positive extremal of
	\begin{equation*}
		\lambda_F(\Omega)=\inf_{\phi\in W^{1,p}_0(\Omega)}\left\{\int_{\Omega}F^p(\nabla \phi)dx:\int_{\Omega}|\phi|^pdx=1\right\}.
	\end{equation*}
   By setting
   \begin{equation}
   	\nonumber
   	\lambda(\Omega;U)=\inf_{\phi\in C_c^1(\Omega)}\left\{\int_{\Omega}\lan D^2H(\nabla U)\nabla \phi,\nabla \phi\ran dx:\int_{\Omega}U^{p-2}|\phi|^2dx=1\right\},
   \end{equation}
    then we have 
    \begin{equation}
    	\nonumber
    	\lambda(\Omega;U)=p(p-1)\lambda_F(\Omega).
    \end{equation}
\end{proposition}
\noindent \emph{Proof.} Firstly, we prove 
\begin{equation}
    \label{ep4}
	\lambda(\Omega;U)\le p(p-1)\lambda_F(\Omega).
\end{equation}
Observe that $U$ solves the equation (\ref{eq27}) with $q=p-1$ and $\mu=\lambda_F(\Omega)$, thus $U\in C^{1,\gamma}(\overline\Omega)$ for some $\gamma\in(0,1)$. Denote $\Omega_\tau=\{x\in\Omega:\text{dist}(x,\partial\Omega) \le\tau\}$. For every $\delta>0$, define the cut-off function $\eta_\delta\in C_c^{\infty}(\Omega)$ such that
\begin{equation}
\eta_\delta=0\quad\text{in }\Omega_{\delta/2},\quad \eta_\delta=1\quad\text{in }\Omega\backslash\Omega_\delta,\quad 0\le\eta_\delta\le1\quad\text{and}\quad |\nabla\eta_\delta|\le\frac{C}{\delta}.  \nonumber
\end{equation}
Note that $|U|\le C\delta$ in $\Omega_\delta$ when $\delta$ is small. Hence, $|U\nabla\eta_\delta|\le C$ in $\Omega_\delta$. Define $U_\delta=U\eta_\delta\in C_c^1(\Omega)$. Noticing that $|U_\delta|\le|U|$, then we have
\begin{equation}
\label{ep5}
\lim_{\delta \rightarrow 0}\int_{\Omega}U_\delta^pdx=\int_{\Omega}U^pdx=1 \quad\text{ and}\quad \lim_{\delta \rightarrow 0}\int_{\Omega}U^{p-2}|U_{\delta}|^2dx=\int_{\Omega}U^pdx=1.
\end{equation}
By direct calculation, we obtain
\begin{equation}
\nabla U_\delta=\eta_\delta\nabla U+U\nabla\eta_\delta.  \nonumber
\end{equation}
Since the definition of $\eta_\delta$, we get $\nabla U_\delta=\nabla U$ in $\Omega\backslash\Omega_\delta$. Moreover, we have $|\nabla U|\le C$ and $|\nabla U_\delta|\le C$ in $\Omega$. Therefore as $\delta\rightarrow0$, we obtain
\begin{equation}
\bigg|\int_\Omega F^p(\nabla U_\delta)dx-\int_\Omega F^p(\nabla U)dx\bigg|\le\bigg|\int_{\Omega_\delta} F^p(\nabla U_\delta)dx\bigg|+\bigg|\int_{\Omega_\delta} F^p(\nabla U)dx\bigg|\rightarrow0  \nonumber
\end{equation}
and
\begin{equation}
\label{eqnew1}
\int_\Omega |\nabla U|^{p-2}|\nabla U_\delta-\nabla U|\big(|\nabla U_\delta|+|\nabla U|\big)dx \le C\int_{\Omega_\delta} |\nabla U|^{p-2}dx\rightarrow0,
\end{equation}
where we utilize the Hopf's Lemma in the last estimate.

Thus, we actually get
\begin{equation}
	\label{ep6}
	\lim_{\delta\rightarrow 0}\int_{\Omega}\lan D^2H(\nabla U)\nabla U_\delta,\nabla U_\delta \ran dx=\int_{\Omega}\lan D^2H(\nabla U)\nabla U,\nabla U \ran dx.
\end{equation}
Indeed, by using Lemma $\ref{lem4}$, we obtain
\begin{equation}
	\nonumber
	\begin{split}
		\bigg|\int_{\Omega}\lan D^2H(\nabla U)\nabla U_\delta,\nabla U_\delta \ran dx-&\int_{\Omega}\lan D^2H(\nabla U)\nabla U,\nabla U \ran dx\bigg|\\
		\le&d_0\int_{\Omega}|\nabla U|^{p-2}|\nabla U_\delta-\nabla U|\big(|\nabla U_\delta|+|\nabla U|\big)dx.
	\end{split}
\end{equation}
Thus, estimate $(\ref{ep6})$ follows by \eqref{eqnew1}. By the homogeneity of $F$, $(\ref{ep5})$ and $(\ref{ep6})$,  
we get 
\begin{equation}
	\nonumber
	\lambda(\Omega;U)\le p(p-1)\lambda_F(\Omega), \quad \text{as} \quad \delta \rightarrow 0.
\end{equation}

For the converse inequality, we recall that $U$ solves
\begin{equation}
\label{ep7}
\int_{\Omega}\lan F^{p-1}(\nabla U)\nabla _{\xi}F(\nabla U), \nabla \phi\ran dx=\lambda_F(\Omega)\int_{\Omega}U^{p-1}\phi dx, \quad \text{for every} \quad \phi\in W^{1,p}_0(\Omega).
\end{equation}
Therefore, by the homogeneity of $F$, we can rewrite \eqref{ep7} as
\begin{equation}
	\label{ep8}
	\int_{\Omega}\lan D^2H(\nabla U)\nabla U, \nabla \phi\ran dx=p(p-1)\lambda_F(\Omega)\int_{\Omega}U^{p-1}\phi dx, \quad \text{for every} \quad \phi\in W^{1,p}_0(\Omega). 
\end{equation}
Next, for any $\epsilon>0$, we select a $\phi_\epsilon\in C_c^1(\Omega)$ such that
\begin{equation}
	\nonumber
	\int_{\Omega}\lan D^2H(\nabla U)\nabla \phi_\epsilon, \nabla \phi_\epsilon\ran dx<\lambda(\Omega;U)+\epsilon \quad \text{and} \quad \int_{\Omega}U^{p-2}\phi_\epsilon^2dx=1.
\end{equation}
Then we test $\phi_\epsilon ^2/U$ on the equation $(\ref{ep8})$, and obtain
\begin{equation}
	\label{ep9}
	\int_{\Omega}\left\lan D^2H(\nabla U)\nabla U, \nabla (\frac{\phi_\epsilon^2}{U})\right\ran dx =p(p-1)\lambda_F(\Omega)\int_{\Omega}U^{p-1}\frac{\phi_\epsilon^2}{U} dx.
\end{equation}
According to the Picone's identity (see \cite[Lemma A.1]{brasco2023uniqueness}), we have
\begin{equation}
	\nonumber
	\begin{aligned}
		\left\lan D^2H(\nabla U)\nabla U, \nabla (\frac{\phi_\epsilon^2}{U})\right\ran=&\left \lan D^2H(\nabla U)\nabla \phi_\epsilon, \nabla \phi_\epsilon \right\ran\\
		&-\left\lan D^2H(\nabla U)\left(\phi_\epsilon \frac{\nabla U}{U}-\nabla \phi_\epsilon\right), \left(\phi_\epsilon\frac{\nabla U}{U}-\nabla \phi_\epsilon\right) \right \ran.
	\end{aligned}
\end{equation}
Inserting this into $(\ref{ep9})$, we get
\begin{equation}
	\nonumber
	\begin{split}
		p(p-1)\lambda_F(\Omega)\int_{\Omega}U^{p-2}\phi_\epsilon^2 dx=&\int_{\Omega}\lan D^2H(\nabla U)\nabla \phi_\epsilon, \nabla \phi_\epsilon \ran dx\\
		&-\int_{\Omega}\left\lan D^2H(\nabla U)\left(\phi_\epsilon\frac{\nabla U}{U}-\nabla \phi_\epsilon\right),\left(\phi_\epsilon\frac{\nabla U}{U}-\nabla \phi_\epsilon\right) \right\ran dx\\
		<&\lambda(\Omega,U)+\epsilon,
	\end{split}
\end{equation}
since $D^2H(\nabla U)$ is positive semi-definite. By our choice of $\phi_\epsilon$ and using the arbitrariness of $\epsilon>0$, we deduce
\begin{equation}
	\nonumber
	p(p-1)\lambda_F(\Omega)\le \lambda(\Omega,U).
\end{equation}
Combining \eqref{ep4} and the last inequality, we finally finish the proof. \hfill $\Box$

\begin{definition}
	For $\frac{3}{2}<p<2$ with the same assumptions in Proposition \ref{pro3.2}, define the weighted Sobolev space
	\begin{equation}
		\nonumber
		X^{1,2}(\Omega;|\nabla U|^{p-2}):=\left\{\phi \in W^{1,2}(\Omega):\int_{\Omega}|\nabla U|^{p-2}|\nabla \phi|^2dx<+\infty\right \},
	\end{equation}
endowed with the natural norm
\begin{equation}
	\nonumber
	||\phi||_{X^{1,2}(\Omega;|\nabla U|^{p-2})}=||\phi||_{L^2(\Omega)}+\left(\int_{\Omega}|\nabla U|^{p-2}|\nabla \phi|^2dx\right)^{\frac{1}{2}}.
\end{equation}
\end{definition}

As discussed, $\int_\Omega |\nabla U|^{p-2}dx$ is bounded when $p>\frac{3}{2}$, thanks to Theorem \ref{thmnew}. Under this situation, we can easily verify that $C_c^1(\Omega)\subset X^{1,2}(\Omega;|\nabla U|^{p-2})$. Accordingly, we define $X^{1,2}_0(\Omega;|\nabla U|^{p-2})$ as the completion of $C_c^1(\Omega)$ with respect to the norm $||\cdot||_{X^{1,2}(\Omega;|\nabla U|^{p-2})}$. Since $U\in C^{1,\gamma}(\overline{\Omega})$, we have
\begin{equation}
	\nonumber
	X^{1,2}_0(\Omega;|\nabla U|^{p-2})\subset W^{1,2}_0(\Omega) \quad\text{for}\quad \frac{3}{2}<p<2,
\end{equation}
with continuous inclusions.

\begin{remark}
\label{rmknew}
Following \cite[Appendix B]{brasco2023uniqueness}, one can similarly verify that
\begin{equation}
X^{1,2}(\Omega;|\nabla U|^{p-2})\cap W_0^{1,2}(\Omega)=X^{1,2}_0(\Omega;|\nabla U|^{p-2}).  \nonumber
\end{equation}
The fact that $|\nabla U|^{p-2}dx$ be a finite positive measure plays a key role in the proof.
\end{remark}

\begin{proposition}
	\label{prop3.5}
	Let $\frac{3}{2}<p<2$ and let $\Omega\subset \mathbb{R}^N$ be an open bounded connected set of class $C^2$. With the same notations of Proposition \ref{pro3.2}, the infimum $\lambda(\Omega;U)$ is uniquely attained on the space  $X^{1,2}_0(\Omega;|\nabla U|^{p-2})$ by the function $U$ or $-U$.
\end{proposition}
\noindent \emph{Proof.} Consider $v\in X^{1,2}_0(\Omega;|\nabla U|^{p-2})$ and $\{v_k\}_{k\in \mathbb{N}} \subset C_c^1(\Omega)$ be such that
\begin{equation}
	\nonumber
	\lim_{k\rightarrow \infty}||v_k-v||_{X^{1,2}(\Omega;|\nabla U|^{p-2})}=0,
\end{equation}
then by Remark \ref{rmknew1} with $q=p-1$, we have
\begin{equation}
\label{epnew}
\lim_{k\rightarrow \infty} \int_{\Omega}U^{p-2}v_k^2 dx=\int_{\Omega}U^{p-2}v^2 dx. 
\end{equation}
Also
\begin{equation}
	\label{ep10}
	\lim_{k\rightarrow \infty}\int_{\Omega}\lan D^2H(\nabla U)\nabla v_k,\nabla v_k\ran dx=\int_{\Omega}\lan D^2H(\nabla U)\nabla v,\nabla v\ran dx.
\end{equation}
Indeed, by Lemma $\ref{lem4}$, it follows
\begin{equation}
	\nonumber
	\big|\lan D^2H(\nabla U)\nabla v_k,\nabla v_k\ran-\lan D^2H(\nabla U)\nabla v,\nabla v\ran\big|\le d_0|\nabla U|^{p-2}|\nabla v_k-\nabla v|(|\nabla v_k|+|\nabla v|).
\end{equation}
Then integrating over $\Omega$ and utilizing H\"older's inequality, we obtain 
\begin{equation}
	\nonumber
	\begin{split}
		\Big|\int_{\Omega}\lan D^2H(\nabla U)&\nabla v_k,\nabla v_k\ran dx-\int_{\Omega}\lan D^2H(\nabla U)\nabla v,\nabla v\ran dx\Big|\\
		\le&\ d_0\int_{\Omega}|\nabla U|^{p-2}|\nabla v_k-\nabla v|(|\nabla v_k|+|\nabla v|)dx\\
		\le&\ d_0\left(\int_\Omega |\nabla U|^{p-2}|\nabla v_k-\nabla v|^2dx\right)^{\frac{1}{2}} \left(\int_\Omega |\nabla U|^{p-2}(|\nabla v_k|+|\nabla v|)^2dx\right)^{\frac{1}{2}}.
	\end{split}
\end{equation}
Noticing that the last term converges to $0$, we derive $(\ref{ep10})$.
Moreover, since $C_c^1(\Omega)$ is dense in $X^{1,2}_0(\Omega;|\nabla U|^{p-2})$, we actually prove by \eqref{epnew} and \eqref{ep10} that
\begin{equation}
	\nonumber
	\lambda(\Omega;U)=\inf_{\phi\in X^{1,2}_0(\Omega;|\nabla U|^{p-2})}\left\{\int_\Omega \lan D^2H(\nabla U)\nabla \phi,\nabla \phi\ran:\int_\Omega U^{p-2}\phi^2dx=1 \right\}.
\end{equation}

In order to show that any minimizer must coincide with either $U$ or $-U$, we suppose on the contrary that there exists another minimizer $v\in X^{1,2}_0(\Omega;|\nabla U|^{p-2})$. By definition, we obtain a sequence $\{v_k\}_{k\in \mathbb{N}}\subset C_c^1(\Omega)$ converging to $v$:
\begin{equation}
	\label{epnew2}
	\lim_{k\rightarrow \infty}\left[\int_\Omega |\nabla U|^{p-2}|\nabla v_k-\nabla v|^2dx+\int_\Omega|v_k-v|^2dx\right]=0.
\end{equation}
Recalling that $U$ satisfies $(\ref{ep8})$, then testing $\phi=v_k^2/U$ on $(\ref{ep8})$ yields 
\begin{equation}
	\label{ep11}
	\begin{split}
		\lambda(\Omega;U)\int_\Omega U^{p-2}v_k^2dx=&\int_\Omega \left\lan D^2H(\nabla U)\nabla U,\nabla \left(\frac{v_k^2}{U}\right) \right\ran dx\\
		=& \int_\Omega \lan D^2H(\nabla U)\nabla v_k,\nabla v_k \ran dx\\
		&-\int_\Omega \left\lan D^2H(\nabla U)\left(v_k\frac{\nabla U}{U}-\nabla v_k\right),\left(v_k\frac{\nabla U}{U}-\nabla v_k\right)\right\ran dx,
	\end{split}
\end{equation}
where the second equality follows again by the Picone's identity. 

In order to pass the limit in \eqref{ep11}, we observe by \eqref{epnew2} that $\{(v_k,\nabla v_k)\}_{k\in \mathbb{N}}$ converges almost everywhere to $(v,\nabla v)$. By applying the Fatou's Lemma, we obtain
\begin{equation}
	\nonumber
	\begin{split}
		\liminf_{k\rightarrow \infty}&\int_\Omega\left\lan D^2H(\nabla U)\left(v_k\frac{\nabla U}{U}-\nabla v_k\right),\left(v_k\frac{\nabla U}{U}-\nabla v_k\right)\right\ran dx\\
		&\ge \int_\Omega\left\lan D^2H(\nabla U)\left(v\frac{\nabla U}{U}-\nabla v\right),\left(v\frac{\nabla U}{U}-\nabla v\right)\right\ran dx.
	\end{split}
\end{equation}
Combining \eqref{epnew}, \eqref{ep10} and the previous inequality, we let $k\to\infty$ in $(\ref{ep11})$ and obtain
\begin{equation}
	\nonumber
	\lambda(\Omega;U)+\int_\Omega \left\lan D^2H(\nabla U)\left(v\frac{\nabla U}{U}-\nabla v\right),\left(v\frac{\nabla U}{U}-\nabla v\right)\right\ran dx \le \lambda(\Omega;U).
\end{equation}
By Assumption F2, $D^2H(\nabla U)$ is positive definite in $\Omega$. Thus, it follows that
\begin{equation*}
	v\frac{\nabla U}{U}-\nabla v=0, \quad \text{a.e. in} \quad \Omega.
\end{equation*}
Since $v\in X^{1,2}_0(\Omega;|\nabla U|^{p-2})\subset W^{1,2}_0(\Omega)$ and $U$ is positive inside $\Omega$, we have $v/U\in W^{1,2}_{loc}(\Omega)$. According to the Leibniz's rule, we have
\begin{equation}
	\nonumber
	\nabla \left(\frac{v}{U}\right)=\frac{1}{U}\left(\nabla v-v\frac{\nabla U}{U}\right)=0, \quad \text{a.e. in} \quad \Omega.
\end{equation}
Since $\Omega$ is connected, this implies that $v/U$ is constant in $\Omega$. Hence, we obtain that $v$ is proportional to $U$. We finally deduce the desired result after the normalization taken.  \hfill $ \Box$

\section{Proof of Theorem \ref{thm2}}
\noindent \emph{Proof of Theorem \ref{thm2}}. We divide the proof into three steps. 

\vspace{0.5em}
\noindent \textbf{Step 1. Linearized equation.} We argue by contradiction: suppose that for $\{q_k\}_{k\in \mathbb{N}}$ with $q_k\searrow p-1$, the equation $(\ref{eq2})$ has (at least) two distinct positive solutions $u_k$ and $v_k$ with $q=q_k$.

Without loss of generality, we assume $||u_k||_{L^\infty(\Omega)}\ge ||v_k||_{L^\infty(\Omega)}$ for each $k\in\mathbb{N}$. Denote $M_k:=||u_k||_{L^\infty(\Omega)}$.
We select a subsequence such that $\lim_{k\rightarrow \infty}{||v_k||_{L^\infty(\Omega)}}/{M_k}=\mu'\in [0,1]$.
Define $\tilde{u}_k={u_k}/{M_k}$, $\tilde{v}_k={v_k}/{M_k}$. Then they satisfy 
\begin{equation}
	\nonumber
	-\Delta^F_p \tilde{u}_k=M^{q_k+1-p}_k\tilde{u}^{q_k}_k \quad \text{and} \quad -\Delta^F_p \tilde{v}_k=M^{q_k+1-p}_k\tilde{v}^{q_k}_k \quad \text{in} \quad \Omega,
\end{equation}
with 
\[\tilde{u}_k=\tilde{v}_k=0 \quad \text{on} \quad \partial \Omega. 
\]
By Lemma $\ref{lem3}$, we can extract a subsequence such that $\tilde{u}_k$ converges to $\tilde{u}$ in $C^{1,\beta}(\overline{\Omega})$, where $\tilde{u}$ is the first eigenfunction of $-\Delta^F_p$ in $\Omega$. By a similar way, we have $\tilde{v}_k$ converges to $\mu'\tilde{u}$ in $C^{1,\beta}(\overline{\Omega})$.
Now, $\tilde{u}_k$ and $\tilde{v}_k$ verify
\begin{equation}
	\label{ee1}
	\int_\Omega \lan  F^{p-1}(\nabla \tilde{u}_k)\nabla F(\nabla \tilde{u}_k),\nabla \phi \ran dx=M^{q_k+1-p}_k\int_\Omega \tilde{u}^{q_k}_k\phi dx,
\end{equation}
and
\begin{equation}
	\label{ee2}
	\int_\Omega \lan  F^{p-1}(\nabla \tilde{v}_k)\nabla F(\nabla \tilde{v}_k),\nabla \phi \ran dx=M^{q_k+1-p}_k\int_\Omega \tilde{v}^{q_k}_k\phi dx,
\end{equation}
for every $\phi \in W^{1,p}_0(\Omega)$. For every $z,w\in \mathbb{R}^N$, we have 
\begin{equation}
	\label{ee3}
	\begin{split}
		 F^{p-1}(z)\nabla F(z)-F^{p-1}(w)\nabla F(w)
		 &=\nabla _{\xi}(\frac{H(\xi)}{p})\bigg |_{\xi=z}-\nabla _{\xi}(\frac{H(\xi)}{p})\bigg |_{\xi=w}\\
		 &=\frac{1}{p}\int_{0}^{1}\frac{d}{dt}\big(\nabla H(tz+(1-t)w)\big)dt\\
		&=\frac{1}{p}\left \lan \int_{0}^{1}(D^2 H(tz+(1-t)w))dt,z-w \right \ran.
	\end{split}
\end{equation}
Similarly, for every $a,b\in[0,+\infty)$, we have
\begin{equation}
	\label{ee4}
	\begin{split}
		a^q-b^q&=\int_{0}^{1}\frac{d}{dt}\big(ta+(1-t)b\big)^qdt\\
		&=q\left(\int_{0}^{1}(ta+(1-t)b)^{q-1}dt\right)(a-b).
	\end{split}
\end{equation}
Then subtracting the equations $(\ref{ee1})$ and $(\ref{ee2})$, and using the identities $(\ref{ee3})$ and $(\ref{ee4})$, we obtain 
\begin{equation}
	\label{ee5}
	\begin{split}
		\int_\Omega\lan A_k(x)\nabla (\tilde u_k-\tilde{v}_k),\nabla \phi \ran dx&=
		M_k^{q_k+1-p}q_k\int_\Omega\int_{0}^{1}(t\tilde{u}_k+(1-t)\tilde{v}_k)^{q_k-1}(\tilde{u}_k-\tilde{v}_k)\phi dtdx,
	\end{split}
\end{equation}
where
\begin{equation}
	\nonumber
	A_k(x):=\frac{1}{p} \int_{0}^{1}D^2H(t\nabla \tilde{u}_k+(1-t)\nabla \tilde{v}_k)dt.
\end{equation}

In the following, we consider the normalization of $\tilde u_k-\tilde{v}_k$ by its weighted $L^2$-norm:
\begin{equation}
    \label{eenew111}
	w_k=\frac{\tilde u_k-\tilde{v}_k}{S_k}\in C^1(\overline\Omega)\cap W^{1,p}_0(\Omega),
\end{equation}
where
\begin{equation}
S_k:=\left(\int_\Omega\int_{0}^{1}(t\tilde{u}_k+(1-t)\tilde{v}_k)^{q_k-1}|\tilde{u}_k-\tilde{v}_k|^2 dtdx\right)^{\frac{1}{2}}>0.  \nonumber
\end{equation}
Then by $(\ref{ee5})$, $w_k$ solves
\begin{equation}
	\label{ee6}
	\begin{split}
		\int_\Omega\lan A_k(x)\nabla w_k,\nabla \phi \ran dx&=M_k^{q_k+1-p}q_k\int_\Omega\int_{0}^{1}(t\tilde{u}_k+(1-t)\tilde{v}_k)^{q_k-1}w_k\phi dtdx.
	\end{split}
\end{equation}
In particular, take $\phi=w_k$ and it yields
\begin{equation}
	\label{ee7}
	\begin{split}
		\int_\Omega\lan A_k(x)\nabla w_k,\nabla w_k \ran dx&=M_k^{q_k+1-p}q_k\int_\Omega\int_{0}^{1}(t\tilde{u}_k+(1-t)\tilde{v}_k)^{q_k-1}w_k^2 dtdx.
	\end{split}
\end{equation}

\noindent \textbf{Step 2. Convergence of $w_k$.} By the definition of $S_k$, we have
\begin{equation}
	\label{ee8}
	\int_\Omega\int_{0}^{1}(t\tilde{u}_k+(1-t)\tilde{v}_k)^{q_k-1}w_k^2 dtdx=1,
\end{equation}
and thus
\begin{equation}
	\label{ee8.5}
	\int_\Omega\lan A_k(x)\nabla w_k,\nabla w_k \ran dx\le C.
\end{equation}
By Assumption F2, we compute directly (note that $p<2$)
\begin{equation}
	\label{ee9}
	\begin{split}
		\lan A_k\xi,\xi\ran &=\frac{1}{p}\int_{0}^{1} \left\lan D^2H(t\nabla \tilde{u}_k+(1-t)\nabla \tilde{v}_k)\xi,\xi \right\ran dt\\
		&\ge\frac{c_0}{p}\left(\int_{0}^{1}|t\nabla \tilde{u}_k+(1-t)\nabla \tilde{v}_k|^{p-2}dt\right)|\xi|^2\\
		&\ge \frac{c_0}{p}\big(|\nabla \tilde{u}_k|+|\nabla \tilde{v}_k|\big)^{p-2}|\xi|^2.
    \end{split}
\end{equation}
Inserting this inequality into \eqref{ee8.5} with $\xi=\nabla w_k$, we have
\begin{equation}
	\nonumber
	\tilde{C}\int_\Omega|\nabla w_k|^2dx\le\int_\Omega (|\nabla \tilde{u}_k|+|\nabla \tilde{v}_k|)^{p-2}|\nabla w_k|^2dx\le C ,\quad \text{for every} \quad k\in\mathbb{N}.
\end{equation}
By Corollary $\ref{cor2.2}$, there exists a $w\in L^2(\Omega)\cap W^{1,2}_0(\Omega)$ such that $\{w_k\}_{k\in \mathbb{N}}$ converges strongly in $L^2(\Omega)$ and weakly in $W^{1,2}(\Omega)$ to $w$, up to a subsequence. In addition, we can similarly obtain
\begin{equation}
    \label{eenew}
\lim_{k\rightarrow\infty}\int_\Omega\int_{0}^{1}(t\tilde{u}_k+(1-t)\tilde{v}_k)^{q_k-1}|w_k-w|^2 dtdx=0.
\end{equation}

Set 
\[Z:=\{x\in \Omega :|\nabla\tilde u|=0\} \quad\text{and} \quad Z_{\delta}:=\{x\in \Omega :|\nabla\tilde u|\le\delta \}.\]
By Theorem \ref{thmnew}, we have $|Z|=0$. Then for every $\epsilon>0$, there exists a $\delta>0$ such that $|Z_\delta|\le\epsilon$. By the $C^1(\overline\Omega)$ convergence of $\tilde u_k$ and $\tilde v_k$, it follows
\begin{equation*}
	A_k\rightarrow A_{\infty}:=\frac{1}{p} \int_{0}^{1}D^2H(t\nabla \tilde{u}+(1-t)\mu'\nabla \tilde{u})dt  \quad \text{uniformly in } Z_{\delta}^c:=\Omega\setminus Z_\delta.
\end{equation*}
Now it is enough to pass the limit in $(\ref{ee6})$. We then assume $\phi\in C_c^1(\Omega)$. The convergence of the right-hand-side of (\ref{ee6}) follows from the uniform convergence of $\tilde{u}_k$, $\tilde{v}_k$ and $L^2$-convergence of $w_k$. For the left-hand-side, we have
\begin{equation}
	\nonumber
	\begin{split}
	&\bigg|\int_\Omega \lan A_k(x)\nabla w_k,\nabla \phi\ran dx-\int_\Omega \lan A_{\infty}(x)\nabla w,\nabla \phi\ran dx\bigg|\\
	\le&\bigg|\int_{Z_\delta^c} \big\lan (A_k(x)-A_{\infty}(x))\nabla w_k,\nabla \phi\big\ran dx\bigg| +\bigg|\int_{Z_\delta^c} \big\lan A_{\infty}(x)(\nabla w_k-\nabla w),\nabla \phi\big\ran dx\bigg|\\
	&+\bigg|\int_{Z_\delta} \lan A_k(x)\nabla w_k,\nabla \phi\ran dx\bigg| +\bigg|\int_{Z_\delta} \lan A_{\infty}(x)\nabla w,\nabla \phi\ran dx\bigg|\\
	\le&||A_k-A_{\infty}||_{L^{\infty}(Z_\delta^c)}||\nabla \phi||_{L^{\infty}(\Omega)} \int_{\Omega} |\nabla w_k|dx +\bigg|\int_{Z_\delta^c} \big\lan A_{\infty}(x)(\nabla w_k-\nabla w),\nabla \phi\big\ran dx\bigg|\\
	&+\bigg|\int_{Z_\delta} \lan A_k(x)\nabla w_k,\nabla \phi\ran dx\bigg| +\bigg|\int_{Z_\delta} \lan A_{\infty}(x)\nabla w,\nabla \phi\ran dx\bigg|\\
	:=& I_1+I_2+I_3+I_4.
	\end{split}
\end{equation}
Since $A_k\rightarrow A_{\infty}$ uniformly in $Z_\delta^c$ and $\{w_k\}_{k\in\mathbb{N}}$ weakly converges to $w$ in $W^{1,2}(\Omega)$, we obtain $I_1$ and $I_2$ tend to zero as $k\rightarrow\infty$. For $I_3$, we utilize H\"older's inequality and get
\begin{equation}
I_3\le\bigg|\int_{Z_\delta} \lan A_k(x)\nabla w_k,\nabla w_k\ran dx\bigg|^{\frac12} \bigg|\int_{Z_\delta} \lan A_k(x)\nabla\phi,\nabla \phi\ran dx\bigg|^{\frac12}.  \nonumber
\end{equation}
By (\ref{ee8.5}) and continuity of integral, the last integral can be arbitrarily small when $\delta\ll1$. Thus $I_3$ tends to zero as $\delta\rightarrow0$. A similar estimate on $I_4$ holds if given $\int_\Omega \lan A_\infty\nabla w,\nabla w\ran dx$ bounded, which can be derived by $(\ref{ee10})$ at below.

Therefore by taking $k\rightarrow\infty$ in (\ref{ee6}), we get
\begin{equation}
\int_\Omega \lan A_{\infty}(x)\nabla w,\nabla \phi\ran dx=(p-1)\lambda_F(\Omega)\int_\Omega\int_{0}^{1}(t\tilde u+(1-t)\mu'\tilde u)^{p-2}w\phi dtdx,  \nonumber
\end{equation}
i.e., 
\begin{equation}
	\nonumber
		\int_{\Omega}  \lan D^2H(\nabla \tilde{u})\nabla w,\nabla \phi\ran dx
		=p(p-1)\lambda_F(\Omega)\int_\Omega\tilde u^{p-2}w\phi dx, \quad\text{for every}\quad \phi\in C_c^1(\Omega).
\end{equation}
In order to pass the limit in (\ref{ee7}), we first denote
\[\Omega_n=\{|\nabla w|\le n\}\cap\{|\nabla \tilde u|\ge \frac1n\}.\]
Then for every $n,k\in \mathbb{N}$, we have
\begin{equation}
	\nonumber
	\begin{split}
		\int_\Omega\lan A_k\nabla w_k,\nabla w_k\ran dx&\ge\int_{\Omega_n}\lan A_k\nabla w_k,\nabla w_k\ran dx\\
		&\ge\int_{\Omega_n}\lan A_k\nabla w,\nabla w\ran dx+2\int_{\Omega_n}\lan A_k\nabla w,\nabla w_k-\nabla w\ran dx.
	\end{split}
\end{equation}
Since $A_k\rightarrow A_{\infty}$ uniformly in $\Omega_n$ and $w_k\rightarrow w$ weakly in $W^{1,2}(\Omega)$, it follows
\begin{equation}
	\nonumber
	\lim_{k\rightarrow \infty}\int_{\Omega_n}\lan A_k\nabla w,\nabla w_k-\nabla w\ran dx=0.
\end{equation}
Thus, we have by Fatou's Lemma
\begin{equation}
	\label{eenew1}
		\liminf_{k\rightarrow \infty}\int_\Omega\lan A_k\nabla w_k,\nabla w_k\ran dx\ge\liminf_{k\rightarrow \infty}\int_{\Omega_n}\lan A_k\nabla w,\nabla w\ran dx =\int_{\Omega_n}\lan A_{\infty}\nabla w,\nabla w\ran dx.
\end{equation}
For the right-hand-side in (\ref{ee7}), we first assume $\{w_k\}_{k\in\mathbb{N}}$ converges almost everywhere to $w$, up to a subsequence. Then by (\ref{ee8}) and Fatou's Lemma again,
\begin{equation}
\int_{\Omega}\int_{0}^{1}(t\tilde{u}+(1-t)\mu\tilde{u})^{p-2}w^2 dtdx \le \liminf_{k\rightarrow\infty}\int_\Omega\int_{0}^{1}(t\tilde{u}_k+(1-t)\tilde{v}_k)^{q_k-1}w_k^2 dtdx=1.  \nonumber
\end{equation}
Since $\tilde u_k$, $\tilde v_k$ converge to $\tilde u$, $\mu\tilde u$ in $C^{1,\beta}(\overline{\Omega})$ and by Hopf's Lemma, we have 
\[\frac{|\tilde u|}{2}\le|\tilde u_k|\le2|\tilde u| \quad\text{and}\quad \frac{\mu|\tilde u|}{2}\le|\tilde v_k|\le2\mu|\tilde u|\]
for $k$ sufficiently large. Thus by using (\ref{eenew}) and the dominated convergence theorem, we obtain
\begin{equation}
	\label{eenew2}
	\begin{split}
&\lim_{k\rightarrow\infty}M_k^{q_k+1-p}q_k\int_\Omega\int_{0}^{1}(t\tilde{u}_k+(1-t)\tilde{v}_k)^{q_k-1}w_k^2 dtdx  \\
=&(p-1)\lambda_F(\Omega)\lim_{k\rightarrow\infty}\int_\Omega\int_{0}^{1}(t\tilde{u}_k+(1-t)\tilde{v}_k)^{q_k-1}w^2 dtdx  \\
=&(p-1)\lambda_F(\Omega)\int_\Omega\int_{0}^{1}(t\tilde{u}+(1-t)\mu\tilde{u})^{p-2}w^2 dtdx.
	\end{split}
\end{equation}
Combining (\ref{ee7}), (\ref{eenew1}), (\ref{eenew2}) and taking $n\rightarrow \infty$, we get
\begin{equation}
	\label{ee10}
	\frac{1}{p}\int_\Omega \lan D^2H(\nabla \tilde{u})\nabla w,\nabla w\ran dx\le (p-1)\lambda_F(\Omega)\int_\Omega\tilde{u}^{p-2}w^2dx.
\end{equation}
Recalling that 
\begin{equation}
	\nonumber
	\lan D^2H(z)\xi,\xi\ran \ge c_0|z|^{p-2}|\xi|^2, \quad \text{for every} \quad z,\xi \in \mathbb{R}^N,
\end{equation}
the estimate $(\ref{ee10})$ implies that $w$ belongs to the weighted Sobolev space $X^{1,2}(\Omega;|\nabla \tilde u|^{p-2})$. $w$ is also nontrivial, by virtue of $(\ref{ee8})$ and $(\ref{eenew})$. Furthermore by Remark \ref{rmknew}, we have
\begin{equation}
	\nonumber
	w\in X^{1,2}(\Omega;|\nabla \tilde u|^{p-2})\cap W^{1,2}_0(\Omega)= X^{1,2}_0(\Omega;|\nabla \tilde u|^{p-2}).
\end{equation}

\noindent\textbf{Step 3. Conclusion.} From the fact that $w\in X^{1,2}_0(\Omega;|\nabla \tilde u|^{p-2})$ is nontrivial together with Proposition $\ref{pro3.2}$, Proposition $\ref{prop3.5}$ and $(\ref{ee10})$, we infer that $w$ must be proportional to either $\tilde{u}$ or $-\tilde{u}$. This illustrates that  $w$ does not change sign.

However, by Lemma $\ref{lem 3.8}$, we notice that $\tilde{u}_k-\tilde{v}_k$ must change sign. Denote by $w^\pm_k$ the positive and negative part of $w_k$, respectively. Then we obtain that each
\begin{equation}
	\nonumber
	\Omega^\pm_k:=\{x\in\Omega:w^\pm_k(x)>0\}
\end{equation}
has positive measure. Note that $w^\pm_k$ are feasible in equation $(\ref{ee6})$, which yields
\begin{equation}
	\label{ee11}
	\int_\Omega \lan A_k\nabla w^\pm_k,\nabla w^\pm_k\ran  dx=M_k^{q_k+1-p}q_k\int_\Omega\int_{0}^{1}(t\tilde{u}_k+(1-t)\tilde{v}_k)^{q_k-1}|w_k^\pm|^2 dtdx.
\end{equation}
Thus, we utilize H\"older's inequality, Theorem $\ref{thm6new}$, $(\ref{ee9})$ and $(\ref{ee11})$ to obtain
\begin{equation*}
	\begin{split}
		&\int_\Omega\int_{0}^{1}(t\tilde{u}_k+(1-t)\tilde{v}_k)^{q_k-1}|w_k^\pm|^2 dtdx \\
		\le& \left(\int_\Omega\int_{0}^{1}(t\tilde{u}_k+(1-t)\tilde{v}_k)^{(q_k-1)\sigma/2}|w_k^\pm|^{\sigma} dtdx\right)^{\frac{2}{\sigma}} |\Omega^\pm_k|^{\frac{\sigma-2}{\sigma}}\\
		\le& \frac{1}{\mathcal{T}}|\Omega^\pm_k|^{\frac{\sigma-2}{\sigma}}\int_{\Omega}|\nabla w^\pm_k|^2dx\\
		\le& \frac{p(2K)^{2-p}}{c_0\mathcal{T}}|\Omega^\pm_k|^{\frac{\sigma-2}{\sigma}}\int_{\Omega}\lan A_k \nabla w^\pm_k,\nabla w^\pm_k \ran dx\\
		\le& \tilde C|\Omega^\pm_k|^{\frac{\sigma-2}{\sigma}}\int_\Omega\int_{0}^{1}(t\tilde{u}_k+(1-t)\tilde{v}_k)^{q_k-1}|w_k^\pm|^2 dtdx,
	\end{split}
\end{equation*}
for some exponent $2<\sigma<\sigma_0$, where $K$ denotes the uniform $C^1$-norm for $\tilde{u}_k$ and $\tilde{v}_k$. This implies
\begin{equation}
	\nonumber
	|\Omega_k^\pm|\ge \frac{1}{C'}, \quad \text{for every} \quad k\in \mathbb{N},
\end{equation}
for some constant $C'$ independent of $k$. We finally reach a contradiction to the fact that $w_k$ strongly converges in $L^2(\Omega)$ to a constant-sign function $w$. This finishes the proof.  \hfill $ \Box$

\section{Proof of Theorem \ref{thm1}}
The constant $Q_{0,p}$ is defined in \cite[Section 4]{bidaut2001nonexistence}.
\begin{lemma}
	\label{lem5.1} 
    Let $\Omega\subset \mathbb{R}^N$ be an open bounded connected set of class $C^2$. For $p-1<q<Q_{0,p}$, let $u\in W^{1,p}_0(\Omega)$ be a positive solution of $(\ref{eq2})$. Then there exists $\delta_1=\delta_1(N,\Omega,p,q)>0$ such that $||u||_{L^\infty(\Omega)}\le C$ holds for every $F$ with $||F(\xi)-|\xi|||_{C^{3,\alpha}(S^{N-1})}\le \delta_1$, where the constant $C$ depends only on $N,p,q,\delta_1$ and $\Omega$.
\end{lemma}
\noindent \emph{Proof.} We argue by contradiction. Suppose that there exists a sequence $\{F_k\}_{k\in \mathbb{N}}$ such that $||F_k(\xi)-|\xi|||_{C^{3,\alpha}(S^{N-1})}\rightarrow0$ and accordingly $\{u_k\}_{k\in\mathbb{N}}\in W^{1,p}_0(\Omega)$ be the solution of $(\ref{eq2})$ with $F=F_k$ satisfying $M_k:=||u_k||_{L^\infty(\Omega)}\rightarrow+\infty$. Let $x_k$ be the point in $\Omega$ where the maximum of $u_k$ is achieved. Define 
\[\tilde{u}_k(y):=\frac{u_k(\mu_ky+x_k)}{M_k},
\]
then $\tilde{u}_k$ satisfies $0\leq\tilde{u}_k\leq1$, $\tilde{u}_k(0)=1$ and 
\begin{equation}
    \label{eqf13}
	\left\{
	\begin{aligned}
		-\Delta^{F_k}_{p} \tilde{u}_k&=M_{k}^{q+1-p}\mu^{p}_k\tilde{u}^{q}_k \quad \text{in} \quad \Omega_k,\\
		\tilde{u}_k&=0 \quad \text{on} \quad \partial \Omega_k,
	\end{aligned}
	\right.
\end{equation}
where $\Omega_k=\{y\in \mathbb{R}^N:\mu_ky+x_k\in \Omega\}. $
Let $M_{k}^{q+1-p}\mu^p_k=1$, then $\mu_k\rightarrow0$ as $k\rightarrow\infty$. Equation \eqref{eqf13} turns into
\begin{equation*}
	\left\{
	\begin{aligned}
		-\Delta^{F_k}_p \tilde{u}_k&=\tilde{u}^{q}_k  \quad \text{in} \quad \Omega_k,\\
		\tilde{u}_k&=0 \quad \text{on} \quad \partial \Omega_k.
	\end{aligned}
	\right.
\end{equation*}
Up to a subsequence, two situations may occur:
\[\text{either} \quad \text{dist}(x_k,\partial \Omega)\mu^{-1}_k\rightarrow +\infty \quad \text{or} \quad \text{dist}(x_k,\partial \Omega)\mu^{-1}_k\rightarrow d \geq 0.
\]

\textbf{Case 1:} $\lim_{k\rightarrow\infty}\text{dist}(x_k,\partial \Omega)\mu^{-1}_k\rightarrow +\infty$. Then
for any $r>1$, $B_r(0)\subset\Omega_k$ holds for $k$ large. According to the local $C^{1,\alpha}$ estimates (see \cite[Theorem 1 \& 2]{dibenedetto1983c1+}), there exists a constant $C=C(r)>0$ uniformly for $k$ such that 
\[||\tilde{u}_k||_{C^{1,\alpha}(B_r(0))}\le C.
\]
By selecting a subsequence, $\tilde{u}_k\rightarrow \tilde{u}$ in $C^{1,\beta}_{loc}(\mathbb{R}^N)$ for some $\beta\in (0,\alpha)$.
Moreover, $\tilde{u}(0)=1$, $\tilde{u}\ge 0$ in $\mathbb{R}^N$, and $\tilde{u}$ solves
\[\int_{\mathbb{R}^N} |\nabla \tilde{u}|^{p-2} \nabla \tilde{u}\cdot\nabla \psi dx =\int_{\mathbb{R}^N}\tilde{u}^{q}\psi dx \quad \text{for all} \quad \psi \in C^{\infty}_c(\mathbb{R}^N),
\]
namely, $\tilde{u}$ is a distributional solution of
\begin{equation*}
	-\Delta _p\tilde{u}=\tilde{u}^q \quad \text{in} \quad \mathbb{R}^N.
\end{equation*} 
By Liouville's Theorem (see \cite[Theorem 3.1]{bidaut2001nonexistence}), we must have $\tilde{u}\equiv0$. But this contradicts $\tilde{u}(0)=1$.

\textbf{Case 2:} $\text{dist}(x^k,\partial \Omega)\mu_k^{-1}\rightarrow d\ge 0$. Then $x^k\rightarrow x^0\in \partial \Omega$. By a similar argument as in Lemma $\ref{lem2}$, we infer that $d>0$ and there exists a subsequence of $\{v_k\}$ such that $v_k\rightarrow v$ in $C^{1,\beta}_{loc}(\mathbb{R}^N_+)$. Moreover, $v$ solves
\begin{equation*}
		\left \{
	\begin{aligned}
		-\Delta_p{v}&={v}^{q} \quad  \text{in} \quad \mathbb{R}^N_+,\\
		v&=0 \quad  \text{on} \quad \partial \mathbb{R}^N_+,\\
        v(de_N)&=1.
	\end{aligned}
\right.
\end{equation*}
By Liouville's Theorem (see \cite[Theorem 4.4]{bidaut2001nonexistence}), we must have $v\equiv0$. But this contradicts $v(de_N)=1$. This lemma is proved. \hfill $ \Box$

\vspace{0.5em}
\noindent \emph{Proof of Theorem \ref{thm1}.} Suppose on the contrary that the equation $(\ref{eq2})$ always admits (at least) two distinct positive solutions. We then take a sequence $\{F_k\}_{k\in \mathbb{N}}$ such that 
\[\lim_{k\rightarrow \infty}||F_k(\xi)-|\xi|||_{C^{3,\alpha}(S^{N-1})}=0.
\]
and for every $k\in\mathbb{N}$, there exist two distinct positive solutions $u_k$ and $v_k$ of $(\ref{eq2})$ with $F=F_k$. By Lemma \ref{lem5.1}, Theorem $\ref{thm5}$ and assumptions of Theorem $\ref{thm1}$, $u_k, v_k \rightarrow u \in C^{1,\beta}(\overline{\Omega})$, up to a subsequence. And $u$ is the unique positive solution of
\begin{equation*}
	\left\{
	\begin{aligned}
		-\Delta_p u&=u^q \quad \text{in} \quad  \Omega,\\
		u&=0 \quad \text{on} \quad \p \Omega.
	\end{aligned}
	\right.
\end{equation*} 
Claim: $u\not\equiv0$. Testing $u_k$ on the equation $(\ref{eq2})$ for $F=F_k$, we obtain
\[\int_{\Omega} |F_k(\nabla u_k)|^pdx=\int_{\Omega} u_k^{q+1}dx.
\]
By Assumption F1, we have 
\[\int_{\Omega} |\nabla u_k|^pdx\le a_0^{-p}\int_{\Omega} u_k^{q+1}dx,
\]
where $a_0>0$ is independent of $k$. By the Sobolev's inequality, we get
\[\left(\int_{\Omega}  u_k^{q+1}dx\right)^{\frac{p}{q+1}}\le a_0^{-p}C \int_{\Omega} u_k^{q+1}dx.
\]
Note that $q>p-1$, thus $u_k$ has a uniform positive lower bound for $L^{q+1}$-norm, which implies that $u\not\equiv0$.

The equations solved by $u_k$ and $v_k$ can be written as 
\begin{equation}
	\label{eee1}
	\int_{\Omega}F_k^{p-1}(\nabla u_k)\nabla _{\xi} F(\nabla u_k)\cdot \nabla \phi= \int_\Omega u^q_k\phi dx,
\end{equation} 
and
\begin{equation}
	\label{eee2}
	\int_{\Omega}F_k^{p-1}(\nabla v_k)\nabla _{\xi} F(\nabla v_k)\cdot \nabla \phi= \int_\Omega v^q_k\phi dx,
\end{equation} 
for every $\phi\in W^{1,p}_0(\Omega)$.
Define $w_k=\frac{u_k-v_k}{S_k}$, where
\[S_k:=\left(\int_{\Omega}\int_{0}^{1}(tu_k+(1-t)v_k)^{q-1}|u_k-v_k|^2 dtdx\right)^{\frac{1}{2}}.\]
Then subtracting the equations \eqref{eee1} and \eqref{eee2}, and using the identities $(\ref{ee3})$ and $(\ref{ee4})$, we obtain
\begin{equation}
	\label{eee5}
	\begin{split}
		\int_{\Omega}\lan A_k(x)\nabla w_k,\nabla \phi \ran dx=
		q\int_{\Omega}\int_{0}^{1}(tu_k+(1-t)v_k)^{q-1}w_k\phi dtdx,
	\end{split}
\end{equation}
where
\[A_k(x):=\frac{1}{p}\int_{0}^{1}D^2H_k(t\nabla u_k+(1-t)\nabla v_k)dt.
\] 
Repeating the process in the proof of Theorem $\ref{thm2}$, we have $\{w_k\}_{k\in \mathbb{N}}$ converges to $w$ strongly in $L^2(\Omega)$ and weakly in $W^{1,2}(\Omega)$. In addition, $w$ is a nontrivial function in $X_0^{1,2}(\Omega;|\nabla u|^{p-2})$. Denote
\[A_{\infty}(x):=(p-2)|\nabla u|^{p-4}\nabla u\otimes\nabla u+|\nabla u|^{p-2}I_N.\]
For every $\phi\in C_c^1(\Omega)$, passing the limit in (\ref{eee5}), we obtain
\begin{equation}
	\nonumber
	\int_{\Omega}\lan A_{\infty}\nabla w,\nabla \phi \ran dx=q\int_{\Omega}u^{q-1}w\phi dx,
\end{equation}
which illustrates that $w\in X_0^{1,2}(\Omega;|\nabla u|^{p-2})$ is a nontrivial solution to the linearized equation of (\ref{eq04}). This leads to a contradiction to the non-degeneracy of $u$. \hfill $\Box $

\vspace{1em}
\noindent\textbf{Data Availability.} 
This paper has no associated data.

\vspace{1em}
\noindent\textbf{Conflict Of Interest.}
The authors declare that there is no conflict of interest.

\bibliographystyle{plain}
\bibliography{citeliouville}

\medskip\medskip
{\em Address and E-mail:}
\medskip\medskip

{\em Rongxun He}

{\em School of Mathematical Sciences, Fudan University, Shanghai 200433, China}

{\em rxhe24@m.fudan.edu.cn}

\medskip\medskip

{\em Wei Ke}

{\em School of Mathematics, Southwest Minzu University, Chengdu 610041, China}

{\em wke25@swun.edu.cn}

\end{document}